\documentclass[10pt, a4paper]{article}
\usepackage[
  paperwidth=7.25in,
  paperheight=9.5in,
  bindingoffset=.75in,
  heightrounded,
]{geometry}

\usepackage{lipsum}

\usepackage{amsmath}
\usepackage{amsfonts}
\usepackage{amssymb}
\usepackage{amsmath,authblk,hyperref}
\newcommand{\RNum}[1]{\uppercase\expandafter{\romannumeral #1\relax}}
\newcommand{\quotes}[1]{``#1"}

\date{}
\begin{document}

\title{\textbf{\textsf{ Common Fixed Point Theorems for Three Transformations on a Vector Valued $S$-metric Spaces}}}

\author[1]{Pooja Yadav\thanks{poojayadav.math.rs@igu.ac.in}}
\author[2]{Mamta Kamra\thanks{mkhaneja15@gmail.com}}
\affil[1, 2]{\quotes{Department of Mathematics, IGU, Meerpur-122502, India.}}
\maketitle
\begin{abstract}
\noindent\textbf{Objectives:} To demonstrate common fixed point results on a complete vector S-metric space.\\
\textbf{Methods:} By using weakly compatible and point of coincidence property of self transformation.\\
\textbf{Findings:} Obtained common fixed point results and application in vector S-metric space and also demonstrated some examples to strengthen our results.\\
\textbf{Novelty:} In this manuscript, we extend S-metric space into vector S-metric space where space is linear lattice valued. By utilizing weakly compatible and point of coincidence property, we demonstrate common fixed point results for three transformations on complete vector S-metric space and  also demonstrate an application in vector S-metric space.Compatible.
\end{abstract}


\section{Introduction}
The fixed point theory is noticeable in its simplicity and can be used to solve problems in various areas like integral equations, differential equations, chemical equations, and many more. In literature,  many researchers have established numerous results in this direction. Several generalizations of metric space are available in the literature, some of which are as follows:  \\
The idea of 2-metric spaces was given by G\"{a}hler\cite{Gahler} in 1963. The researchers  in  \cite{Ha} have discovered results in 2-metric space. In 2006,  Mustafa and Sims \cite{Mustafa} gave a more general concept known as G-metric.
In 2012, Sedghi et al. \cite{Sedghi} yielded S-metric space which is extension of metric space. Many researchers established common fixed point(CFP) results in S-metric space with different contractive conditions. P. Yadav and M. Kamra \cite{Pooja} extended the idea of S-metric space in vector metric space and thereby extend the existing results for scalar cases. In this article, we prove CFP theorems for three transformations on  complete vector S-metric spaces(VSMS). 
\section{Definitions and Preliminaries}
We start with some basic mathematical definitions and established findings.\\
\textbf{\emph{Definition 2.1}} \cite{Pooja}
A partially ordered set(poset) is called a lattice if each set with two elements has an infimum and a supremum.\\ 
\textbf{\emph{Example 2.2}} Let $S=\{1, 2\}$ and \[P(S)= \{\emptyset, \{1\}, \{2\}, \{1,2\}\} \]is a lattice, where $\emptyset$ is the minimum element and $\{1, 2\}$ is the maximum element.\\
\textbf{Definition 2.3} \cite{Pooja}
Let $\complement$ be  linear space which is real and $(\complement, \preceq)$ be a poset. Then the poset $(\complement, \preceq)$ is said to be an ordered linear space if it follows the  properties mentioned below:
\begin{itemize}
\item[(a)]$ \Theta_1 \preceq \Theta_2 \Longrightarrow \Theta_1+\Theta_3 \preceq \Theta_2+\Theta_3$
\item[(b)] $\Theta_1 \preceq \Theta_2 \Longrightarrow \omega \Theta_1 \preceq \omega \Theta _2$\hspace{1cm}$\forall \Theta_1, \Theta_2, \Theta_3 \in \complement$ and $\omega >0$
\end{itemize}  
\textbf{\emph{Definition 2.4}} \cite{Pooja}
An ordered linear space which is a lattice is called a linear lattice. This is also called a Riesz space.\\ \\
\textbf{\emph{Example 2.5}} The space of all continuous functions $C[0, 1]$ is a linear lattice.\\ 
\textbf{\emph{Lemma 2.6}} If $V$ is a linear lattice and $ \alpha\preceq\gamma\alpha$ where $\alpha \in V^+$ and $0\leq\gamma <1 $, then $\alpha=0$\\\\
\textbf{\emph{Definition 2.7}} \cite{Pooja}
A linear lattice V is said to be Archimedean if $\dfrac{1}{n}\xi\downarrow 0$ for each $\xi\in V^+$ where\[V^+= \{\xi \in V:\xi \succeq 0\}\]
\textbf{\emph{Definition 2.8}} \cite{Pooja}
Let $V$ be a linear lattice and $\Re$ be a nonvoid set. A transformation $S:\Re\times \Re\times \Re\rightarrow V$ on $\Re$ is called vector $S$-metric if it satisfies the following conditions:
\begin{itemize}
\item [(a)] $S(\varrho_1, \varrho_2, \varrho_3) \succeq 0$,
\item[(b)]$S(\varrho_1, \varrho_2, \varrho_3) = 0$ iff $\varrho_1 =\varrho_2 = \varrho_3$,
\item[(c)]$S(\varrho_1, \varrho_2, \varrho_3) \preceq S(\varrho_1, \varrho_1, \beta)+S(\varrho_2, \varrho_2, \beta)+S(\varrho_3, \varrho_3, \beta)$,
\end{itemize}
\,\,for all  $\varrho_1, \varrho_2, \varrho_3, \beta \in \Re$.\\
Then $(\Re, S, V)$ is called vector $S$-metric space(VSMS).\\\\
\textbf{Example 2.09} Let $V$ be a linear lattice and $\Re$ be a nonvoid set. A transformation $S:\Re\times \Re\times \Re\rightarrow V$ is defines as \[S(\varrho_1, \varrho_2, \varrho_3)= max\{S(\varrho_1, \varrho_1, \varrho_2), S(\varrho_2, \varrho_2, \varrho_3), S(\varrho_3, \varrho_3, \varrho_1)\}\] $\forall \varrho_1, \varrho_2, \varrho_3\in \Re$. Then $(\Re, S, V)$ is a VSMS.\\ \\
\textbf{Definition 2.10} \cite{Pooja}
In a vector $S$-metric space $(\Re, S, V)$,  $\langle \vartheta_n \rangle$ is called $V$-convergent to some $\vartheta \in V$ if there is a sequence $\langle \mu_n \rangle$ in $V$ satisfying $\mu_n \downarrow 0$ and $S(\vartheta_n, \vartheta_n, \vartheta) \leq \mu_n$ and denote it by $\vartheta_n \xrightarrow{S,V} \vartheta$.\\ \\
\textbf{Definition 2.11}\cite{Pooja}
In a vector $S$-metric space $(\Re, S, V)$,  $\langle \vartheta_n \rangle$ is known as $V$-Cauchy sequence if there is a sequence $\langle \mu_n \rangle$ in $V$ satisfying $\mu_n \downarrow 0$ and $S(\vartheta_n, \vartheta_n, \vartheta_{n+q}) \leq \mu_n$ holds for all $q$ and $n$.\\\\
\textbf{Definition 2.12} \cite{Pooja}
A VSMS is called $V$-complete if each $V$-Cauchy sequence in $\Re$ is $V$-convergent to a limit in $\Re$.\\ \\
\textbf{\emph{Definition 2.13}} \cite{Shahzad}
Let $p$ and $q$ be two transformations in a vector $S$-metric space $(\Re, S, V)$. Then the transformations are called weakly compatible (WC) if $pq\alpha=qp\alpha$ whenever $p\alpha=q\alpha$ for $\alpha \in \Re$. \\ \\
\textbf{\emph{Definition 2.14}} \cite{Ali}
Let $p, q:\Re\rightarrow \Re$ be two transformations. Then $p$ and $q$ are said to have a coincidence point(CP ) $\alpha \in \Re$ if $p(\alpha)=q(\alpha)=\beta$ and $\beta$ is called a point of coincidence(PC) of $p$ and $q$. \\ \\
\textbf{Lemma 2.15} \cite{Harish} Let $\Re$ be a set, $p$ and $q$ are WC transformations of $\Re$. If $\omega=p\alpha=q\alpha$ for some $\alpha\in \Re$ where $p$ and $q$ have $\omega$ as PC which is unique then $\omega$ is CFP of $p$ and $q$ which is unique.\\
\textbf{\emph{Lemma 2.16}} \cite{Rajpal}
If $(\Re, S, V)$ is a VSMS, then  \[S(\vartheta, \vartheta, \mu)=S(\mu, \mu, \vartheta) \,\,\,\,\,\,\,\forall \mu, \vartheta \in \Re.\]
\section{Main Results}
With the existing research work for scalar valued cases, we pursue new generalized CFP results for three transformations on V-complete VSMS.\\
\textbf{Theorem 3.1} Let $(\Re, S, V)$ be V-complete VSMS and $V$-Archimedean. Assume that $p, q, k:\Re\rightarrow \Re$ satisfy the conditions stated below:
\begin{itemize}
\item[(i)] $\forall \xi, \gamma\in \Re$
\begin{eqnarray}
S(p\xi, p\xi, q\gamma) &\preceq& \hslash_1 S(k\xi, k\xi, k\gamma)+ \hslash_2 S(p\xi, p\xi, k\xi)+ \hslash_3 S(q\gamma, q\gamma, k\gamma)\nonumber\\&& +\hslash_4 S(p\xi, p\xi, k\gamma)+\hslash_5 S(q\gamma, q\gamma, k\xi)
\end{eqnarray}
where $\hslash_i$ for i=1, 2, \dots, 5 are positive constant with $2\hslash_1+2\hslash_2+2\hslash_3+4\hslash_4+4\hslash_5<1;$
\item[(ii)] $p(\Re)\cup q(\Re)\subset k(\Re)$;
\item[(iii)] one of $p(\Re)$, $q(\Re)$ or $k(\Re)$ is a V-complete subspace of $\Re$.
\end{itemize} 
Then $\{p, k\}$ and $\{q, k\}$ have unique PC in $\Re$. If  $\{p, k\}$ and $\{q, k\}$ are WC, then $p$, $q$ and $k$ have a CFP in $\Re$ which is unique. \\
\textbf{Proof.} Suppose $\xi_0\in\Re$. Since $p(\Re)\subset k(\Re)$, $\exists\,\,\, \xi_1 \in \Re$ such that $p(\xi_0)=k(\xi_1)=\gamma_1$ and  $q(\Re)\subset k(\Re)$, $\exists\,\,\, \xi_2 \in \Re$ such that $q(\xi_1)=k(\xi_2)=\gamma_2$. \\
Continue in this manner\[\exists\,\,\,\, \xi_{2b+1}\in \Re \,\,\,\, s.t. \,\,\,\, \gamma_{2b+1}=p\xi_{2b}=k\xi_{2b+1}\]\[\,\,\,\,\exists\,\,\,\,\, \xi_{2b+2}\in \Re \,\,\,\, s.t. \,\,\,\, \gamma_{2b+2}=q\xi_{2b+1}=k\xi_{2b+2}\] for b=0, 1, \dots.
Firstly, we assert that \\
\begin{eqnarray}\label{eq3}
S(\gamma_{2b+1}, \gamma_{2b+1}, \gamma_{2b+2}) \preceq \alpha S(\gamma_{2b}, \gamma_{2b}, \gamma_{2b+1})
\end{eqnarray} 
$\forall$ b. We have\\
\begin{eqnarray}
S(\gamma_{2b+1}, \gamma_{2b+1}, \gamma_{2b+2})&=& S(p\xi_{2b}, p\xi_{2b}, q\xi_{2b+1})\nonumber\\
&\preceq& \hslash_1 S(k\xi_{2b}, k\xi_{2b}, k\xi_{2b+1})+ \hslash_2 S(p\xi_{2b}, p\xi_{2b}, k\xi_{2b})+ \nonumber\\&&\hslash_3 S(q\gamma_{2b+1}, q\gamma_{2b+1}, k\gamma_{2b+1}) +\hslash_4 S(p\xi_{2b}, p\xi_{2b}, k\gamma_{2b+1})\nonumber\\&&+\hslash_5 S(q\gamma_{2b+1}, q\gamma_{2b+1}, k\xi_{2b})\nonumber
\end{eqnarray}
\begin{eqnarray}
&\preceq& \hslash_1 S(\gamma_{2b}, \gamma_{2b}, \gamma_{2b+1})+ \hslash_2 S(\gamma_{2b+1}, \gamma_{2b+1}, \gamma_{2b})+ \nonumber\\&&\hslash_3 S(\gamma_{2b+2}, \gamma_{2b+2}, \gamma_{2b+1}) +\hslash_4 S(\gamma_{2b+1}, \gamma_{2b+1}, \gamma_{2b+1})\nonumber\\&&+\hslash_5 S(\gamma_{2b+2}, \gamma_{2b+2}, \gamma_{2b})\nonumber\\
S(\gamma_{2b+1}, \gamma_{2b+1}, \gamma_{2b+2})&\preceq& \hslash_1 S(\gamma_{2b}, \gamma_{2b}, \gamma_{2b+1})+ \hslash_2 S(\gamma_{2b}, \gamma_{2b}, \gamma_{2b+1})+ \nonumber\\&&\hslash_3 S(\gamma_{2b+1}, \gamma_{2b+1}, \gamma_{2b+2})+\hslash_5 [2S(\gamma_{2b+2}, \gamma_{2b+2}, \gamma_{2b+1})\nonumber\\&&+S(\gamma_{2n}, \gamma_{2n}, \gamma_{2n+1})]\nonumber
\end{eqnarray}
\begin{eqnarray}
(1-\hslash_3-2\hslash_5)S(\gamma_{2b+1}, \gamma_{2b+1}, \gamma_{2b+2})&\preceq& (\hslash_1+\hslash_2+ \hslash_5)S(\gamma_{2b}, \gamma_{2b}, \gamma_{2b+1})\nonumber\\
S(\gamma_{2b+1}, \gamma_{2b+1}, \gamma_{2b+2})&\preceq&\dfrac{\hslash_1+\hslash_2+\hslash_5}{1-\hslash_3-2\hslash_5}S(\gamma_{2b}, \gamma_{2b}, \gamma_{2b+1})\nonumber
\end{eqnarray}
So, 
\begin{eqnarray}\label{eq7}
S(\gamma_{2b+1}, \gamma_{2b+1}, \gamma_{2b+2}) \preceq \alpha S(\gamma_{2b}, \gamma_{2b}, \gamma_{2b+1})
\end{eqnarray}
 where $\alpha=\dfrac{\hslash_1+\hslash_2+\hslash_5}{1-\hslash_3-2\hslash_5}\prec 1.$ \\
Similiarly
\begin{eqnarray}
S(\gamma_{2b+3}, \gamma_{2b+3}, \gamma_{2b+2})&=&S(p\xi_{2b+2}, p\xi_{2b+2}, q\xi_{2b+1})\nonumber\\
&\preceq& \hslash_1 S(k\xi_{2b+2}, k\xi_{2b+2}, k\xi_{2b+1})+\hslash_2 S(p\xi_{2b+2}, p\xi_{2b+2}, k\xi_{2b+2})\nonumber\\&&+\hslash_3 S(q\xi_{2b+1}, q\xi_{2b+1})+\hslash_4S(p\xi_{2b+2}, p\xi_{2b+2}, k\xi_{2b+1})\nonumber\\&&+\hslash_5 S(q\xi_{2b+1}, q\xi_{2b+1}, k\xi_{2b+2})\nonumber\\
&=& \hslash_1 S(\gamma_{2b+2}, \gamma_{2b+2}, \gamma_{2b+1})+\hslash_2S(\gamma_{2b+3}, \gamma_{2b+3}, \gamma_{2b+2})+ \nonumber\\&& \hslash_3 S(\gamma_{2b+2}, \gamma_{2b+2}, \gamma_{2b+1})+\hslash_4 S(\gamma_{2b+3}, \gamma_{2b+3}, \gamma_{2b+1})+\nonumber\\&&\hslash_5 S(\gamma_{2n+2}, \gamma_{2n+2}, \gamma_{2n+2})\nonumber\\
&=& \hslash_1 S(\gamma_{2b+2}, \gamma_{2b+2}, \gamma_{2b+1})+\hslash_2S(\gamma_{2b+3}, \gamma_{2b+3}, \gamma_{2b+2})+ \nonumber\\&& \hslash_3 S(\gamma_{2b+2}, \gamma_{2b+2}, \gamma_{2b+1})+\hslash_4[2S(\gamma_{2n+3}, \gamma_{2n+3}, \gamma_{2n+2})+\nonumber\\&& S(\gamma_{2b+1}, \gamma_{2b+1}, \gamma_{2b+2})]\nonumber
\end{eqnarray}
\begin{eqnarray} 
&=& \hslash_1 S(\gamma_{2b+2}, \gamma_{2b+2}, \gamma_{2b+1})+\hslash_2S(\gamma_{2b+3}, \gamma_{2b+3}, \gamma_{2b+2})+ \nonumber\\&& \hslash_3 S(\gamma_{2b+2}, \gamma_{2b+2}, \gamma_{2b+1})+\hslash_4[2S(\gamma_{2b+3}, \gamma_{2b+3}, \gamma_{2b+2})+\nonumber\\&& S(\gamma_{2b+2}, \gamma_{2b+2}, \gamma_{2b+1})]\nonumber
\\
(1-\hslash_2-2\hslash_4)S(\gamma_{2b+3}, \gamma_{2b+3}, \gamma_{2b+2})&\preceq& (\hslash_1+\hslash_3+\hslash_4)S(\gamma_{2b+2}, \gamma_{2b+2}, \gamma_{2n+1})\nonumber\\
S(\gamma_{2b+3}, \gamma_{2b+3}, \gamma_{2b+2})&\preceq& \dfrac{\hslash_1+\hslash_3+\hslash_4}{1-\hslash_2-2\hslash_4}S(\gamma_{2b+2}, \gamma_{2b+2}, \gamma_{2n+1})\nonumber
\end{eqnarray}
Thus, 
\begin{eqnarray}\label{eq8}
S(\gamma_{2b+3}, \gamma_{2b+3}, \gamma_{2b+2}) \preceq \alpha S(\gamma_{2b+2}, \gamma_{2b+2}, \gamma_{2n+1})
\end{eqnarray}
 where \,\,\,\,$\alpha=\dfrac{\hslash_1+\hslash_3+\hslash_4}{1-\hslash_2-2\hslash_4}\prec1.$ \\ 
 From \eqref{eq7} and \eqref{eq8}, we have
 \[S(\gamma_{b}, \gamma_{b}, \gamma_{b+1})\preceq \alpha^b S(\gamma_0, \gamma_0, \gamma_1)\]
 It can be proved same as in Theorem 3.1, we get that $\langle \gamma_b\rangle$ is $V$-Cauchy sequence. Assume that $k(\Re)$ is V-complete. Then $\exists\,\,\,\,\,\gamma\in k(\Re)$ such that \[k\xi_{2b}=\gamma_{2b}\xrightarrow{S, V}\gamma\,\,\,\,\,\, and \,\,\,\,k\xi_{2b+1}=\gamma_{2b+1}\xrightarrow{S, V}\gamma.\]
 Thus $\exists\,\,\,\, \langle c_b\rangle\in V$ such that $c_b\downarrow0$ and \[S(k\xi_{2b}, k\xi_{2b}, \gamma)\preceq c_b \,\,\,\,and\,\,\,\,\,  S(k\xi_{2b+1}, k\xi_{2b+1}, \gamma)\preceq c_{b+1}.\]
 Since $k:\Re\rightarrow\Re$ be a transformation on $\Re$, $\exists$\,\,\,\, $\omega \in\Re$ such that $k\omega=\gamma$.\\
 Now we assert that $p\omega=\gamma.$ For this, consider
 \begin{eqnarray}
 S(p\omega, p\omega, \gamma)&\preceq& 2S(p\omega, p\omega,q\xi_{2b+1})+S(\gamma, \gamma, q\xi_{2b+1})\nonumber\\
&=& 2S(p\omega, p\omega,q\xi_{2b+1})+S(q\xi_{2b+1}, q\xi_{2b+1}, \gamma)\nonumber\\ 
&\preceq& 2[\hslash_1 S(k\omega, k\omega, k\xi_{2b+1})+ \hslash_2 S(p\omega, p\omega, k\omega)+\nonumber\\&&\hslash_3 S(q\xi_{2b+1}, q\xi_{2b+1}, k\xi_{2b+1})+\hslash_4S(p\omega, p\omega, k\xi_{2b+1})+\nonumber\\&& \hslash_5 S(q\xi_{2b+1}, q\xi_{2b+1}, k\omega)]+S(q\xi_{2b+1}, q\xi_{2b+1}, \gamma)\nonumber\\
&\preceq& 2\hslash_1 S(\gamma, \gamma, k\xi_{2b+1})+ 2\hslash_2 S(p\omega, p\omega, \gamma)+\nonumber\\&&2\hslash_3 [S(q\xi_{2b+1}, q\xi_{2b+1}, \gamma)+S(k\xi_{2b+1}, k\xi_{2b+1}, \gamma)]+\nonumber\\&&2\hslash_4 [S(p\omega, p\omega, \gamma)+S(k\xi_{2b+1}, k\xi_{2b+1}, \gamma)]\nonumber\\&& +2\hslash_5 S(q\xi_{2b+1}, q\xi_{2b+1}, \gamma)+S(q\xi_{2b+1}, q\xi_{2b+1}, \gamma)\nonumber
\\
&=& 2\hslash_1 S(k\xi_{2b+1}, k\xi_{2b+1}, \gamma)+ 2\hslash_2 S(p\omega, p\omega, \gamma)+\nonumber
\end{eqnarray}
\begin{eqnarray} 
&&2\hslash_3 [S(q\xi_{2b+1}, q\xi_{2b+1}, \gamma) +S(k\xi_{2b+1}, k\xi_{2b+1}, \gamma)]+\nonumber\\&&2\hslash_4 [S(p\omega, p\omega, \gamma)+S(k\xi_{2b+1}, k\xi_{2b+1}, \gamma)] +\nonumber\\&&2\hslash_5 S(q\xi_{2b+1}, q\xi_{2b+1}, \gamma)+S(q\xi_{2b+1}, q\xi_{2b+1}, \gamma)\nonumber
\\ 
&=& (2\hslash_1+2\hslash_3+2\hslash_4)S(k\xi_{2b+1}, k\xi_{2b+1}, \gamma)+\nonumber\\&&(4\hslash_3+2\hslash_5+1)S(q\xi_{2b+1}, q\xi_{2b+1}, \gamma)\nonumber\\&&+(2\hslash_2+4\hslash_4)S(p\omega, p\omega, \gamma)\nonumber
\\
(1-2\hslash_2-4\hslash_4) S(p\omega, p\omega, \gamma)&\preceq&(2\hslash_1+2\hslash_3+2\hslash_4)S(k\xi_{2b+1}, k\xi_{2b+1}, \gamma)+\nonumber\\&&(4\hslash_3+2\hslash_5+1)S(q\xi_{2b+1}, q\xi_{2b+1}, \gamma)\nonumber\\
&\preceq&(2\hslash_1+2\hslash_3+2\hslash_4)S(k\xi_{2b+1}, k\xi_{2b+1}, \gamma)+\nonumber\\&&(4\hslash_3+2\hslash_5+1)S(k\xi_{2b+2}, k\xi_{2b+2}, \gamma)\nonumber
\\
&\preceq&(2\hslash_1+2\hslash_3+2\hslash_4)c_{b+1}+(4\hslash_3+2\hslash_5+1)c_{b+2}\nonumber
\\
&\preceq& (2\hslash_1+2\hslash_3+2\hslash_4)c_{b}+(4\hslash_3+2\hslash_5+1)c_{b}\nonumber\\
S(p\omega, p\omega, \gamma)&\preceq& \dfrac{2\hslash_1+6\hslash_3+2\hslash_4+2\hslash_5+1}{1-2\hslash_2-4\hslash_4}c_{b}\nonumber
 \end{eqnarray}
 for all $b$. So, $S(p\omega, p\omega, \gamma)=0$, i.e. $p\omega=\gamma$. Thus $p\omega=k\omega=\gamma$, i.e. $\gamma$ is PC of the transformations $p$ and $k$, $\omega$ is CP of the transformations $p$ and $k$. 
\begin{eqnarray}
S(q\omega, q\omega, \gamma) &\preceq& 2S(q\omega, q\omega, p\xi_{2b})+S(\gamma, \gamma, p\xi_{2b}) \nonumber
\\
&=&2S(p\xi_{2b}, p\xi_{2b}, q\omega)+S(p\xi_{2b}, p\xi_{2b}, \gamma) \nonumber\\ 
&\preceq& 2\hslash_1 S(k\xi_{2b}, k\xi_{2b}, k\omega)+2\hslash_2 S(p\xi_{2b}, p\xi_{2b}, k\xi_{2b})+\nonumber\\&&2\hslash_3 S(q\omega, q\omega, k\omega)+2\hslash_4 S(p\xi_{2b}, p\xi_{2b}, k\omega)+\nonumber\\&&2\hslash_5 S(q\omega, q\omega, k\xi_{2b})+S(k\xi_{2b+1}, k\xi_{2b+1}, \gamma)\nonumber
\\
&\preceq& 2\hslash_1 S(k\xi_{2b}, k\xi_{2b}, \gamma)+2\hslash_2 [2S(p\xi_{2b}, p\xi_{2b}, \gamma)+\nonumber\\&&S(k\xi_{2b}, k\xi_{2b}, \gamma)]+2\hslash_3 S(q\omega, q\omega, \gamma)+\nonumber\\&&2\hslash_4 S(p\xi_{2b}, p\xi_{2b}, \gamma)+2\hslash_5[ 2S(q\omega, q\omega, \gamma)\nonumber\\&&+S(k\xi_{2b}, k\xi_{2b}, \gamma)]+S(k\xi_{2b+1}, k\xi_{2b+1}, \gamma)\nonumber
\\
&=& 2\hslash_1 S(k\xi_{2b}, k\xi_{2b}, \gamma)+2\hslash_2 [2S(k\xi_{2b+1}, k\xi_{2b+1}, \gamma)+\nonumber
\\
&&S(k\xi_{2b}, k\xi_{2b}, \gamma)]+2\hslash_3 S(q\omega, q\omega, \gamma)+\nonumber\\&&2\hslash_4 S(k\xi_{2b+1}, k\xi_{2b+1}, \gamma)+2\hslash_5[ 2S(q\omega, q\omega, \gamma)+\nonumber\\&&S(k\xi_{2b}, k\xi_{2b}, \gamma)]+S(k\xi_{2b+1}, k\xi_{2b+1}, \gamma)\nonumber\\
&\preceq& (1+4\hslash_2+2\hslash_4)S(k\xi_{2b+1}, k\xi_{2b+1}, \gamma)+\nonumber\\&&(2\hslash_1+2\hslash_2+2\hslash_5)S(k\xi_{2b}, k\xi_{2b}, \gamma)+\nonumber\\&&(2\hslash_3+4\hslash_5)S(q\omega, q\omega, \gamma)\nonumber
\end{eqnarray}
\begin{eqnarray} 
(1-2\hslash_3-4\hslash_5)S(q\omega, q\omega, \gamma) &\preceq&  (1+4\hslash_2+2\hslash_4)c_{b+1}+(2\hslash_1+2\hslash_2+2\hslash_5)c_b\nonumber\\
&\preceq&  (1+4\hslash_2+2\hslash_4)c_{b}+(2\hslash_1+2\hslash_2+2\hslash_5)c_b\nonumber\\
&=& (1+2\hslash_1+6\hslash_2+2\hslash_4+2\hslash_5)c_b\nonumber\\
S(q\omega, q\omega, \gamma) &\preceq& \dfrac{1+2\hslash_1+6\hslash_2+2\hslash_4+2\hslash_5}{1-2\hslash_3-4\hslash_5}c_b\nonumber
\end{eqnarray} 
for all b. So, $S(q\omega, q\omega, \gamma)=0$ i.e. $q\omega=\gamma$.\\
Thus $q\omega=k\omega=\gamma$, $\gamma$ is PC of the transformations $q$ and $k$, $\omega$ is CP of transformations $q$ and $T$. Now, we shall assert that $\gamma$ is unique PC of the pairs $\{p, k\}$ and $\{q, k\}$. Let $ \Theta'$ be another PC of these three transformations, then 
\[p\omega'= q\omega'=k\omega'=\Theta'\,\,\, for \,\,\, \omega'\in \Re.\] 
Now, we have 
\begin{eqnarray}
S(\gamma, \gamma, \Theta')&=& S(p\omega, p\omega, q\omega')\nonumber\\
&\preceq& \hslash_1 S(k\omega, k\omega, k\omega')+\hslash_2S(p\omega, p\omega, k\omega)+\hslash_3 S(q\omega', q\omega', k\omega')+ \nonumber\\&& \hslash_4S(p\omega, p\omega, k\omega')+\hslash_5S(q\omega', q\omega', k\omega)\nonumber\\
&=& \hslash_1 S(\gamma, \gamma, \Theta')+\hslash_2S(\gamma, \gamma, \gamma)+\hslash_3 S(\gamma', \gamma', \gamma')+ \nonumber\\&& \hslash_4S(\gamma, \gamma, \gamma')+\hslash_5S(\Theta', \Theta', \gamma)\nonumber
\\
S(\gamma, \gamma, \Theta')&\preceq& (\hslash_1+\hslash_4+\hslash_5)S(\gamma, \gamma, \Theta')\nonumber
\end{eqnarray}
Hence $S(\gamma, \gamma, \Theta')=0$, i.e. $\gamma=\Theta'$.\\
If $\{p, k\}$ and $\{q, k\}$ are WC, then $\gamma$ is CFP of $p$, $q$, and $k$ by lemma 2.15. The proof for the rest cases are similiar when $q(\Re)$ and $k(\Re)$ are complete. \\\\
We get the following corollary which is deduction of the above Theorem. In this corollary, we demonstrate CFP result in VSMS by utilize two transformations instead of three transformations .\\
\textbf{Corollary 3.2} Let $(\Re, S, V)$ be V-complete VSMS and $V$-Archimedean. Assume that $p, k:\Re\rightarrow \Re$ satisfy the conditions stated below:
\begin{itemize}
\item[(i)] $\forall \xi, \gamma\in \Re$
\begin{eqnarray}
S(p\xi, p\xi, p\gamma) &\preceq& \hslash_1 S(k\xi, k\xi, k\gamma)+ \hslash_2 S(p\xi, p\xi, k\xi)+ \hslash_3 S(p\gamma, p\gamma, k\gamma)\nonumber\\&& +\hslash_4 S(p\xi, p\xi, k\gamma)+\hslash_5 S(p\gamma, p\gamma, k\xi)
\end{eqnarray}
where $\hslash_i$ for i=1, 2, \dots, 5 are positive constant with $2\hslash_1+2\hslash_2+2\hslash_3+4\hslash_4+4\hslash_5<1;$
\item[(ii)] $p(\Re)\subset k(\Re)$;
\item[(iii)] one of $p(\Re)$ or $k(\Re)$ is a V-complete subspace of $\Re$.
\end{itemize} 
Then $\{p, k\}$  has unique PC in $\Re$. If  $\{p, k\}$ is WC, then $p$ and $k$ have a CFP in $\Re$ which is unique. 
\section{Application}
We present an application of corollary 3.2 . To do this, we establish the CFP theorem on a V-complete VSMS.\\
\textbf{Theorem 4.1} Let $(\Re, S, V)$ be V-complete VSMS and $V$-Archimedean. Assume that $p, k:\Re\rightarrow \Re$ satisfy the conditions stated below:
\begin{itemize}
\item[(i)] 
\begin{eqnarray} \label{eq11}
\int_{0}^{S(p\xi, p\xi, p\gamma)}\mho (\ell)d\ell
 &\preceq& \hslash_1 \int\limits_0^{S(k\xi, k\xi, k\gamma)}\mho (\ell)d\ell+ \hslash_2 \int\limits_0^{S(p\xi, p\xi, k\xi)}\mho (\ell)d\ell\nonumber
\\&&+ \hslash_3 \int\limits_0^{S(p\gamma, p\gamma, k\gamma)}\mho (\ell)d\ell +\hslash_4 \int\limits_0^{S(p\xi, p\xi, k\gamma)}\mho (\ell)d\ell\nonumber\\&&+\hslash_5 \int\limits_0^{S(p\gamma, p\gamma, k\xi)}\mho (\ell)d\ell
\end{eqnarray}
$\forall \xi, \gamma\in \Re$ and $\hslash_1, \dots ,\hslash_5\geq 0$ with $ 2\hslash_1+2\hslash_2+2\hslash_3+4\hslash_4+4\hslash_5<1$ and $\mho:\mathbb{R^+}\cup 0\rightarrow \mathbb{R^+}\cup 0$ is a Lebesgue integrable transformation which is summable , positive and $\int\limits_0^{\epsilon}\mho (\ell)d\ell>0$ for each $\epsilon>0$;
\item[(ii)] $p(\Re)\subset k(\Re)$;
\item[(iii)] p or k is continuous;
\end{itemize}
 Then the transformations  $p$ and $k$ have CFP which is unique in $\Re$ provided p and k are WC.\\
\textbf{Proof.} Let $\xi_0 \in \Re$ and we can take $\gamma_b=p\xi_b=k\xi_{b+1}$, $b=0,1,2,\dots$ \\
From (\ref{eq11}), we have
\begin{eqnarray} 
\int\limits_0^{S(p\xi_b, p\xi_b, p\xi_{b+1})}\mho (\ell)d\ell
 &\preceq& \hslash_1 \int\limits_0^{S(k\xi_b, k\xi_b, k\xi_{b+1})}\mho (\ell)d\ell+ \hslash_2 \int\limits_0^{S(p\xi_b, p\xi_b, k\xi_b)}\mho (\ell)d\ell\nonumber\\&&+ \hslash_3 \int\limits_0^{S(p\xi_{b+1}, p\xi_{b+1}, k\xi_{b+1})}\mho (\ell)d\ell +\hslash_4 \int\limits_0^{S(p\xi_b, p\xi_b, k\xi_{b+1})}\mho (\ell)d\ell\nonumber
 \end{eqnarray}
\begin{eqnarray}
 &&+\hslash_5 \int\limits_0^{S(p\xi_{b+1}, p\xi_{b+1}, k\xi_b)}\mho (\ell)d\ell \nonumber
\\
  &\preceq&\hslash_1 \int\limits_0^{S(\gamma_{b-1}, \gamma_{b-1}, \gamma_{b})}\mho (\ell)d\ell+ \hslash_2 \int\limits_0^{S(\gamma_{b}, \gamma_{b}, \gamma_{b-1})}\mho (\ell)d\ell\nonumber
\\
  &&+ \hslash_3 \int\limits_0^{S(\gamma_{b+1}, \gamma_{b+1}, \gamma_{b})}\mho (\ell)d\ell +\hslash_4 \int\limits_0^{S(\gamma_{b}, \gamma_{b}, \gamma_{b})}\mho (\ell)d\ell\nonumber
\\
&&+\hslash_5 \int\limits_0^{S(\gamma_{b+1}, \gamma_{b+1}, \gamma_{b-1})}\mho (\ell)d\ell \nonumber\\ 
  \int\limits_0^{S(p\xi_b, p\xi_b, p\xi_{b+1})}\mho (\ell)d\ell&\preceq& (\hslash_1+\hslash_2) \int\limits_0^{S(\gamma_{b}, \gamma_{b}, \gamma_{b-1})}\mho (\ell)d\ell+\hslash_3 \int\limits_0^{S(\gamma_{b+1}, \gamma_{b+1}, \gamma_{b})}\mho (\ell)d\ell \nonumber\\&&+\hslash_5 \int\limits_0^{S(\gamma_{b+1}, \gamma_{b+1}, \gamma_{b-1})}\mho (\ell)d\ell \nonumber
\end{eqnarray}
Since
\begin{eqnarray}\label{eq12}
S(\gamma_{b+1}, \gamma_{b+1}, \gamma_{b-1})\preceq 2S(\gamma_{b+1}, \gamma_{b+1}, \gamma_b)+ S(\gamma_{b-1}, \gamma_{b-1}, \gamma_{b})\nonumber\\
S(\gamma_{b+1}, \gamma_{b+1}, \gamma_{b-1})\preceq 2S(\gamma_{b+1}, \gamma_{b+1}, \gamma_b)+ S(\gamma_{b}, \gamma_{b}, \gamma_{b-1}) 
\end{eqnarray}
By using eq(\ref{eq12}) and lemma 2.16, we get
\begin{eqnarray}\label{eq13}
\int\limits_0^{S(p\xi_b, p\xi_b, p\xi_{b+1})}\mho (\ell)d\ell&\preceq& (\hslash_1+\hslash_2) \int\limits_0^{S(\gamma_{b}, \gamma_{b}, \gamma_{b-1})}\mho (\ell)d\ell+\hslash_3 \int\limits_0^{S(\gamma_{b+1}, \gamma_{b+1}, \gamma_{b})}\mho (\ell)d\ell \nonumber
\end{eqnarray}
\begin{eqnarray} 
&&+\hslash_5 \int\limits_0^{2S(\gamma_{b+1}, \gamma_{b+1}, \gamma_b)}\mho (\ell)d\ell +\hslash_5 \int\limits_0^{ S(\gamma_{b}, \gamma_{b}, \gamma_{b-1})} \mho (\ell)d\ell \nonumber
\\ 
\int\limits_0^{S(\gamma_{b}, \gamma_{b}, \gamma_{b+1})}\mho (\ell)d\ell
 &\preceq& (\hslash_1+\hslash_2+\hslash_5) \int\limits_0^{S(\gamma_{b}, \gamma_{b}, \gamma_{b-1})}\mho (\ell)d\ell \nonumber
\\
 && +(\hslash_3+2\hslash_5) \int\limits_0^{S(\gamma_{b}, \gamma_{b}, \gamma_{b+1})}\mho (\ell)d\ell \nonumber
\\
 \int\limits_0^{S(\gamma_{b}, \gamma_{b}, \gamma_{b+1})}\mho (\ell)d\ell
 &\preceq& (\hslash_1+\hslash_2+\hslash_5) \int\limits_0^{S(\gamma_{b-1}, \gamma_{b-1}, \gamma_{b})}\mho (\ell)d\ell \nonumber
\\
 && +(\hslash_3+2\hslash_5) \int\limits_0^{S(\gamma_{b}, \gamma_{b}, \gamma_{b+1})}\mho (\ell)d\ell \nonumber
\\
 \Big(1-(\hslash_3+2\hslash_5)\Big)\int\limits_0^{S(\gamma_{b}, \gamma_{b}, \gamma_{b+1})}\mho (\ell)d\ell
 &\preceq& (\hslash_1+\hslash_2+\hslash_5) \int\limits_0^{S(\gamma_{b-1}, \gamma_{b-1}, \gamma_{b})}\mho (\ell)d\ell \nonumber\\
 \int\limits_0^{S(\gamma_{b}, \gamma_{b}, \gamma_{b+1})}\mho (\ell)d\ell
 &\preceq&\dfrac{\hslash_1+\hslash_2+\hslash_5}{1-(\hslash_3+2\hslash_5)}\int\limits_0^{S(\gamma_{b-1}, \gamma_{b-1}, \gamma_{b})}\mho (\ell)d\ell \nonumber
\\
  \int\limits_0^{S(p\xi_b, p\xi_b, p\xi_{b+1})}\mho (\ell)d\ell&\preceq& \dfrac{\hslash_1+\hslash_2+\hslash_5}{1-(\hslash_3+2\hslash_5)}\int\limits_0^{S(p\xi_{b-1}, p\xi_{b-1}, p\xi_{b})}
  \mho (\ell)d\ell \nonumber\\
   \int\limits_0^{S(p\xi_b, p\xi_b, p\xi_{b+1})}\mho (\ell)d\ell&\preceq& \vartheta \int\limits_0^{S(p\xi_{b-1}, p\xi_{b-1}, p\xi_{b})}
  \mho (\ell)d\ell
\end{eqnarray}
where $\vartheta=\dfrac{\hslash_1+\hslash_2+\hslash_5}{1-(\hslash_3+2\hslash_5)}<1.$\\
From inequality (\ref{eq13}), we get
\[ \int\limits_0^{S(p\xi_b, p\xi_b, p\xi_{b+1})}\mho (\ell)d\ell \preceq \vartheta^b \int\limits_0^{S(p\xi_{0}, p\xi_{0}, p\xi_{1})}\mho (\ell)d\ell. \]
For all $b, m \in \mathbb{N}$, $b<m$, we have
\begin{eqnarray}
  \int\limits_0^{S(\gamma_b, \gamma_b, \gamma_m)}\mho (\ell)d\ell &\preceq& \int\limits_0^{2S(\gamma_b, \gamma_b, \gamma_{b+1})}\mho (\ell)d\ell+\int\limits_0^{S(\gamma_m, \gamma_m, \gamma_{b+1})}\mho (\ell)d\ell \nonumber
\\ 
&\preceq& \int\limits_0^{2S(\gamma_b, \gamma_b, \gamma_{b+1})}\mho (\ell)d\ell+\int\limits_0^{S(\gamma_{b+1}, \gamma_{b+1}, \gamma_m)}\mho (\ell)d\ell \nonumber
\end{eqnarray}
\begin{eqnarray}
  &\preceq& \int\limits_0^{2S(\gamma_b, \gamma_b, \gamma_{b+1})}\mho (\ell)d\ell+\int\limits_0^{2S(\gamma_{b+1}, \gamma_{b+1}, \gamma_{b+2})}\mho (\ell)d\ell + \nonumber\\&&\int\limits_0^{S(\gamma_{b+2}, \gamma_{b+2}, \gamma_m)}\mho (\ell)d\ell \nonumber
\\ 
&&  \vdots\nonumber
\\
&\preceq& \int\limits_0^{2S(\gamma_b, \gamma_b, \gamma_{b+1})}\mho (\ell)d\ell+\int\limits_0^{2S(\gamma_{b+1}, \gamma_{b+1}, \gamma_{b+2})}\mho (\ell)d\ell +\nonumber\\&&\dots+ \int\limits_0^{S(\gamma_{m-1}, \gamma_{m-1}, \gamma_m)}\mho (\ell)d\ell \nonumber
\\
&<& \int\limits_0^{2S(\gamma_b, \gamma_b, \gamma_{b+1})}\mho (\ell)d\ell+\int\limits_0^{2S(\gamma_{b+1}, \gamma_{b+1}, \gamma_{b+2})}\mho (\ell)d\ell +\nonumber\\&&\dots+ \int\limits_0^{2S(\gamma_{m-1}, \gamma_{m-1}, \gamma_m)}\mho (\ell)d\ell \nonumber\\
&<& 2(\vartheta^{b}+\vartheta^{b+1}+\dots+ \vartheta^{m-1})  \int\limits_0^{S(\gamma_0, \gamma_0, \gamma_1)}\mho (\ell)d\ell \nonumber\\
&<& \dfrac{2\vartheta^{b}}{1-\vartheta}  \int\limits_0^{S(\gamma_0, \gamma_0, \gamma_1)} \mho (\ell)d\ell\rightarrow 0 \,\,\,\,as\,\,\,\, b, m\rightarrow \infty \nonumber
\end{eqnarray}
Thus \[\lim_{b, m\rightarrow \infty } S(\gamma_b, \gamma_b, \gamma_m)=0.\]
So, $\langle \gamma_b\rangle \in \Re$. Since $(\Re, S, V)$ is V-complete VSMS, so $\exists \,\,\,f \in \Re$ such that \[\lim_{b \rightarrow \infty} \gamma_b=\lim_{b \rightarrow \infty} p\xi_b= \lim_{b \rightarrow \infty} k\xi_{b+1}=f.\]
Since the transformations $p$ or $k$ is continuous, so we assume that $k$ is continuous, then \[\lim_{b \rightarrow \infty} kp\xi_b= \lim_{b \rightarrow \infty} kk\xi_{b+1}=kf.\] 
Also, $p$ and $k$ are WC, therefore $$\lim_{b \rightarrow \infty} S(pk\xi_b, pk\xi_b, kp\xi_b)=0,$$ this implies \[\lim_{b \rightarrow \infty}pk\xi_b=kf. \]   
From (\ref{eq11}), we have
\begin{eqnarray}
 \int\limits_0^{S(pk\xi_b, pk\xi_b, p\xi_b)}\mho (\ell)d\ell &\preceq&\hslash_1 \int\limits_0^{S(kk\xi_b, kk\xi_b, k\xi_{b})}\mho (\ell)d\ell+ \hslash_2 \int\limits_0^{S(pk\xi_b, pk\xi_b, kk\xi_b)}\mho (\ell)d\ell\nonumber\\&&+ \hslash_3 \int\limits_0^{S(p\xi_{b}, p\xi_{b}, k\xi_{b})}\mho (\ell)d\ell +\hslash_4 \int\limits_0^{S(pk\xi_b, pk\xi_b, k\xi_{b})}\mho (\ell)d\ell\nonumber\\&&+\hslash_5 \int\limits_0^{S(p\xi_{b}, p\xi_{b}, kk\xi_b)}\mho (\ell)d\ell. \nonumber
\end{eqnarray} 
This gives $kf=f$ when $b \rightarrow\infty.$\\
Again from (\ref{eq11}), we have 
\begin{eqnarray} 
\int\limits_0^{S(p\xi_b, p\xi_b, pf)}\mho (\ell)d\ell
 &\preceq& \hslash_1 \int\limits_0^{S(k\xi_b, k\xi_b, kf)}\mho (\ell)d\ell+ \hslash_2 \int\limits_0^{S(p\xi_b, p\xi_b, k\xi_b)}\mho (\ell)d\ell\nonumber\\&&+ \hslash_3 \int\limits_0^{S(pf, pf, kf)}\mho (\ell)d\ell +\hslash_4 \int\limits_0^{S(p\xi_b, p\xi_b, kf)}\mho (\ell)d\ell\nonumber\\&&+\hslash_5 \int\limits_0^{S(pf, pf, k\xi_b)}\mho (\ell)d\ell \nonumber
    \end{eqnarray}
  This gives $pf=f$ when $b \rightarrow\infty.$ So $kf=pf=f$. Thus $f$ is a CFP  of $p$ and $k$.\\
  Now we assert that $f$ is the unique CFP of $p$ and $k$. Let  $\sigma$ be another CFP of the transformations $p$ and $k$, then 
  \begin{eqnarray} 
\int\limits_0^{S(f, f, \sigma)}\mho (\ell)d\ell
 &=& \int\limits_0^{S(pf, pf, p\sigma)}\mho (\ell)d\ell\nonumber
 \end{eqnarray}
\begin{eqnarray}
  &\preceq& \hslash_1 \int\limits_0^{S(kf, kf, k\sigma)}\mho (\ell)d\ell+ \hslash_2 \int\limits_0^{S(pf, pf, kf)}\mho (\ell)d\ell\nonumber
\\
  &&+ \hslash_3 \int\limits_0^{S(p\sigma, p\sigma, k\sigma)}\mho (\ell)d\ell +\hslash_4 \int\limits_0^{S( pf, pf, k\sigma)}\mho (\ell)d\ell\nonumber  
\\    
  &&+\hslash_5 \int\limits_0^{S(p\sigma, p\sigma, kf)}\mho (\ell)d\ell \nonumber
\\ 
 &\preceq& (\hslash_1+\hslash_4) \int\limits_0^{S(f, f, \sigma)}\mho (\ell)d\ell+\hslash_5 \int\limits_0^{S(\sigma, \sigma, f)}\mho (\ell)d\ell\nonumber
\\
  &\preceq& (\hslash_1+\hslash_4+\hslash_5) \int\limits_0^{S(f, f, \sigma)}\mho (\ell)d\ell.\nonumber
    \end{eqnarray}
    Since $\hslash_1+\hslash_4+\hslash_5<1$. So, $S(f, f, \sigma)=0$. Thus $f=\sigma$.\\
    Now, we give an e.g. to demonstrate Theorem 4.1.\\
    \textbf{Example 4.2} Suppose $0\leq\Re\leq 1 $ and a VSM \,\,\,\,$S$ on $\Re\times\Re\times\Re$ be defined as \[S(\alpha, \mu, \varphi)= \vert \alpha-\mu\vert+\vert \mu-\varphi\vert+ \vert \varphi-\alpha \vert\,\,\,\, \alpha, \mu, \varphi\in \Re.\] Then $(\Re, S, V)$ be VSMS. \\We define $p\xi=\dfrac{\xi}{12}$ and $k\xi= \dfrac{\xi}{3}$ and the transformation $p$ is continuous and $p(\Re)\subseteq k(\Re)$. Also 
    \[\int\limits_0^{S(p\xi, p\xi, p\gamma)}\mho (\ell)d\ell\leq \vartheta \int\limits_0^{S(k\xi, k\xi, k\gamma)}\mho (\ell)d\ell\,\,\,\,\,\forall \,\,\,\xi, \gamma\in \Re, \,\,\,\, \vartheta\in \Big[\dfrac{1}{3}, 1\Big] \] and 0 is unique CFP of $p$ and $k$.
\section{Conclusion}
In this manuscript, we demonstrate CFP theorems for three transformations in V-complete VSMS by using WC and PC. Our results are the generalization of existing results of \cite{Rad} from vector metric spaces to vector S-metric spaces. These results shall motivate researchers to solve problems in various areas like differential equations, functional analysis and many more.
\section*{Declaration of Conflicts Interest}
There are no conflicts of interest with this work disclosed by the authors.
\section*{Acknowledgements}
The first author gratefully acknowledges financial support from the "University Grant Commission(UGC), New Delhi, India", under student i'd- 1159/ (CSIRNETJUNE2019).
\section*{Author Contribution Statement}
Pooja Yadav: Identify and defining the problem's, analyzing and then solving.
Mamta Kamra:  Provide guidance and supervising the proposed work for accuracy and clarity.
Pooja Yadav and Mamta Kamra wrote this manuscript.
\section*{Data Availability Statement}
There are no related data with this manuscript.

\end{document}